\documentclass[journal]{IEEEtran}
\usepackage{amsmath}
\usepackage{amsthm}
\usepackage{amsfonts}
\usepackage{bbm} 

\newtheorem{lem}{Lemma}

\newcommand{\IM}{\mathrm{Im}}
\newcommand{\Proj}{\mathrm{Proj}}

\usepackage[color=yellow, textsize=small, textwidth=\marginparwidth, draft]{todonotes}   

\begin{document}

\title{On Certain Properties of Convex Functions}
\author{Miel Sharf and Daniel Zelazo
\thanks{M. Sharf and D. Zelazo are with the Faculty of Aerospace Engineering, Israel Institute of Technology, Haifa, Israel.
    {\tt\small msharf@tx.technion.ac.il, dzelazo@technion.ac.il}}
}

\maketitle
\begin{abstract}
This note deals with certain properties of convex functions.  We provide results on the convexity of the set of minima of these functions, the behaviour of their subgradient set under restriction, and optimization of these functions over an affine subspace.
\end{abstract}

\section{Introduction}
This paper deals with certain properties of convex function, most of them well-known but generally unmentioned in the literature. This work employs subgradient calculus, calculating the subgradient space after restricting to a subspace, or optimizing on a ``moving" affine subspace, and also deals with the collection of minima of the convex function.
\section{The Lemmas}
This paper contains the proof of the three following lemmas:
\begin{lem} \label{Lemma1}
Let $f:\mathbb{R}^n \rightarrow \mathbb{R}$ be a convex function, and let $S:\mathbb{R}^n \rightarrow\mathbb{R}^d$ be some linear operator. Fix some $\zeta \in Im(S)$ and define a map $g:\{y:\: Sy=\zeta\}\to \mathbb{R}$ by $g(x)=f(x)$. Then
\begin{enumerate}
\item $g$ is a convex function on the smaller space $\IM(S^T)$,
\item the subdifferential of $g$ at $x$ is given by
\[
\partial g(x) = \Proj_{\IM(S)} (\partial f(x)),
\]
where $\Proj_W$ is the orthogonal projection on the subspace $W$.
\end{enumerate}
\end{lem}

\begin{lem} \label{Lemma2}
Let $f:\mathbb{R}^d \rightarrow \mathbb{R}$ be a convex function, and let $S:\mathbb{R}^n\rightarrow\mathbb{R}^d$ be some linear operator. Define a map $h:\IM(S^T)\to \mathbb{R}$ by
\[
h(x) = \min_{r: S^Tr=x} f(r),
\]
assuming that the minimum is always achieved. Then
\begin{enumerate}
\item $h$ is a convex function,
\item if $f$ is strictly convex, then $h$ is strictly convex.
\end{enumerate}
\end{lem}

\begin{lem} \label{Lemma3}
Let $C\subseteq \mathbb{R}^n$ be a convex set and let $f:C\rightarrow\mathbb{R}$ be convex. Suppose that $f$ achieves its minimum $m$ in $C$, and let $M=\{x:\: f(x)=m\}$ be the set of $f$'s minima. Then $M$ is convex.
\end{lem}

\section*{Proofs}

We start by proving lemma \ref{Lemma2}:
\begin{proof}
We denote $V = \IM(S^T)$ for simplicity. We take some $x,y\in V$ and $t\in [0,1]$. Our goal is to show that
\[
h(tx+(1-t)y) \le th(x) + (1-t)h(y)
\]
pick $r_x,r_y\in \mathbb{R}^d$ such that $g(x)=f(r_x)$ and $g(y)=f(r_y)$ (these exist by assumption that the minimum is always achieved). Then on one hand, we have $S^T(tr_x+(1-t)r_y) = tx+(1-t)y$ by linearity, so $h(tx+(1-t)y) \le f(tr_x+(1-t)r_y)$. On the other hand, by convexity:
\begin{equation} \label{Lemma2ConvexityInequality}
f(tr_x+(1-t)r_y) \le tf(r_x) + (1-t)f(r_y) = th(x) + (1-t)h(y)
\end{equation}
so we get the wanted inequality by chaining the two inequalities.

As for strict convexity, we should note that if $x\neq y$ then $r_x\neq r_y$, so the inequality \ref{Lemma2ConvexityInequality} becomes strict. This completes the proof of the lemma.
\end{proof}

We now prove lemma \ref{Lemma3}:
\begin{proof}
Let $x,y\in M$ and let $t\in[0,1]$. We need to show that $tx+(1-t)y\in M$. Indeed, because $f$ is convex,
\[
f(tx+(1-t)y) \le tf(x) + (1-t)f(y) = tm+(1-t)m = m
\]
but on the other hand, $f(tx+(1-t)y)$ cannot be smaller than $m$, as $m$ is the minimum of $f$. Thus $f(tx+(1-t)y) = m$ and thus $tx+(1-t)y\in C$.
\end{proof}

Lastly, we prove lemma \ref{Lemma1}, which is the "toughest" of the three:
\begin{proof}
We denote $X=\{y:\: Sy=\zeta\}$, $V=\ker{S}$ and $U=\IM(S^T)$. We know that $\IM(S^T)=\ker(S)^\perp$, so we can identify $\mathbb{R}^n$ as a direct sum of $U$ and $V$. Thus we get a function $F:U\times V\rightarrow\mathbb{R}$ defined by $F(u,v)=f(u+v)$.

In \cite{Rockafeller}, one shows that if $\chi:\mathbb{R}\rightarrow\mathbb{R}$ is convex, then $\partial\chi(t) = [\chi^{\prime}_-,(t)\chi^{\prime}_+(t)]$, where $\chi^{\prime}_{\pm}$ are the one-sided derivatives. Furthermore, we know (again by \cite{Rockafeller}) that if $\rho:\mathbb{R}^n\rightarrow\mathbb{R}$ is convex, and we fix $v_0,v_1\in V$ and define $\chi(t)= \rho(v_0+tv_1)$, then $\chi$'s one-sided derivatives are given via:
\begin{align}
\chi^{\prime}_+(t) &= \max_{\alpha\in\partial\rho(v_0+tv_1)} \alpha^Tv_1 \\ \nonumber
\chi^{\prime}_-(t) &= \min_{\alpha\in\partial\rho(v_0+tv_1)} \alpha^Tv_1 \nonumber
\end{align}
These facts, together with the fact that $\partial\rho(v_0+tv_1)$ is convex, imply that $\partial\chi(t) = v^T\cdot\partial\rho(v_0+tv_1)$

Now, we can finally begin out proof. By assumption, there's some $y\in \mathbb{R}$ such that $Sy=\zeta$. We can decompose $y$ as $u_0+v_0$ for some $u_0\in \IM(S^T)$ and $v_0\in\ker(S)$. The set $X$ is equal to $u_0+V$. Thus, we map $g$ can be described as $g(v)=F(u_0,v)$.

Take some $v_0,v_1\in V$. Restricting $g$ to the line $\{v_0+tv_1: t\in\mathbb{R}\}$ is identical to restricing $f$ to the line $\{(u_0,v_0+tv_1): t\in\mathbb{R}\}$. Thus they yield the same subdifferential sets at $t=0$. By above, we get that:
\[
v_1^T \cdot\partial g(v_0) = v_1^T \cdot \partial F(u_0,v_0) = v_1^T \cdot \Proj_V(\partial F(u_0,v_0))
\]
meaning that the sets $\partial g(v_0)$ and  $\Proj_V(\partial F(u_0,v_0))$ look the same when hit by a linear functional on $V$. However, both of these sets are both convex and closed (see \cite{Rockafeller}). Thus, the separating hyperplanen theorem (see \cite{Conway}) implies that they are equal. Reading what $V$ and $F$ are, we get that for any $x\in \IM(S)$,
\[
\partial g(x) = \Proj_{\IM(S)} (\partial f(x))
\]
which completes the proof.
\end{proof}

\bibliographystyle{ieeetr}
\bibliography{Appendix}

\end{document}